\documentstyle{amsppt}
\tolerance=1000
\global\newcount\headno\global\newcount\subheadno
\global\newcount\eqnno\global\newcount\itemno

\define\showheads{\the\headno.}
\define\itemnum{\showheads\the\itemno}
\define\eqnnum{\showheads\the\eqnno}
\define\newhead{\global\advance\headno by1 \the\headno
  \global\subheadno=0\global\itemno=0\global\eqnno=0}
\define\newsubhead{\global\advance\subheadno by1
  \the\headno.\the\subheadno}

\define\newitem{\global\advance\itemno by1 \itemnum}
\define\neweqn{\global\advance\eqnno by1 \eqnnum}
\define\xname#1#2{\xdef#1{#2}}
\define\itemname#1{\xname#1\itemnum} \let\name\itemname
\define\eqnname#1{\xname#1\eqnnum}

\topmatter
\author Eleanor G. Rieffel \endauthor
\title Groups Quasi-isometric to $\bold H^2\times\bold R$ \endtitle
\address FX Palo Alto Laboratory, 3400 Hillview Ave., Bldg 4, Palo Alto, CA 94304, USA \endaddress
\email rieffel\@pal.xerox.com \endemail
\thanks This research was partially supported by National Science Foundation
grant DMS90-14482.
\endthanks
\endtopmatter

\mag=1200
\pagewidth{6.4 truein}
\pageheight{9 truein}

\define\h2r{\bold H^2\times\bold R}
\define\USL2r{\widetilde{SL}(2,\bold R)}
\define\uth2{UT\bold H^2}
\define\uuth2{\widetilde{UT}\bold H^2}
	$$\text{In memory of Herbert Busemann.}$$

\head Introduction
\endhead


The most powerful geometric tools are those of differential
geometry, but to apply such techniques to finitely generated groups 
seems hopeless at first glance since the natural metric on a finitely
generated group is discrete. However Gromov recognized 
that a group can metrically resemble a manifold in such a 
way that geometric results about that manifold carry over to the
group [18, 20]. This resemblance is formalized in the concept of a 
``quasi-isometry.'' This paper contributes to an ongoing program to
understand which groups are quasi-isometric to which simply connected,
homogeneous, Riemannian manifolds [15, 18, 20] by proving 
that any group quasi-isometric to $\h2r$ is a finite extension of a 
cocompact lattice in $Isom(\h2r)$ or $Isom(\USL2r)$.
\proclaim{Theorem}For any group $\Gamma$ quasi-isometric to the hyperbolic plane
cross the real line, there is an exact sequence
	$$0 @>>> A @>>> \Gamma @>>> G @>>> 0$$
where $A$ is virtually infinite cyclic and $G$ is a finite
extension of a cocompact Fuchsian group.
\endproclaim

With this result, the question of which
finitely generated groups are quasi-isometric to each of Thurston's 
eight geometries [33][31] remains open only for Sol. 
Tukia [35] for $n>2$ and Gabai [14], or
alternatively Casson--Jungreis [6], and Tukia [36] in 
dimension 2 show that any group 
quasi-isometric to $\bold H^n$ can be realized as a finite 
extension of a discrete
cocompact subgroup of the isometries of $\bold H^n$.
The analogous result for $\h2r$ is not true,
because $\h2r$ is quasi-isometric to $\USL2r$, the universal cover
of $SL(2,\bold R)$. 
Gromov's
work on groups of polynomial growth [17] implies that any group
quasi-isometric to $\bold R^n$ is a finite extension of 
a discrete, cocompact
subgroup of the isometry group of $\bold R^n$. 
Similarly for Nil. 
Kapovich and Leeb [22] show that quasi-isometries preserve the
geometric decompostion of Haken manifolds. Their paper, which was
written after the present paper, contains another proof, along
quite different lines, of the main result in this paper.
In their 1996 preprint [24], Kleiner and Leeb generalized the 
main results of the present paper to products of simply connected
nilpotent Lie group with a symmetric space of non-compact type with
no Euclidean de Rham factors.
Chow [5] has shown that every group 
quasi-isometric to the complex hyperbolic plane is a finite
extension of a discrete, cocompact subgroup of the isometry group
of the complex hyperbolic plane.
Pansu [27] has shown that in
quaternionic and Cayley hyperbolic space every quasi-isometry
is within bounded distance of an isometry, so that in fact every group
quasi-isometric to one of these spaces  is a finite extension of a
group naturally isomorphic to a discrete cocompact subgroup of the
isometry group. 
A series of papers by Pansu [27], Schwartz [29, 30],
Farb-Schwartz [13], Kleiner-Leeb [23], Eskin-Farb [8], and Eskin [7]
succeeded in classifying lattices in semi-simple Lie groups up to
quasi-isometry. See [10] for a survey. 
Of course, the question may be asked for
metric spaces other than symmetric spaces or even manifolds. Papers by 
Gromov and Stallings give the result for trees [19][32]. Farb and Mosher
[11][12] prove qusi-isometric rigidity for the solvable Baumslag-Solitar
groups.

The heart of the paper is contained is Sections 3 and 5.
Section 1 gives some background and definitions and sets the notation
for the paper.
Section 2 describes a ``quasi-action" of a group $\Gamma$
quasi-isometric to $\h2r$ on $\h2r$.
Section 3 uses geometric arguments to obtain an action of $\Gamma$
on $\partial\bold H^2$ by quasisymmetric maps.
Section 4 describes the application of results of 
Hinkkanen [21], Gabai [14], Casson--Jungreis [6], 
and Tukia [36] to show that this action is conjugate by a quasisymmetric 
map to an action by a M\" obius group. Section 5 uses largely algebraic
methods to show that the resulting action on $\bold H^2$ by 
isometries must have been properly discontinuous, and thus there is a map 
$\Phi:\Gamma\to G$ where $G$ is a discrete, cocompact 
subgroup of the isometries of $\bold H^2$.  Section 6 
uses geometric arguments from Section 3 to show that $\ker\Phi$ 
must be quasi-isometric to $\bold R$.

The aim of Section 3 is to show that the image 
of a horizontal hyperbolic plane in $\h2r$ under any quasi-isometry induced by
$\Gamma$ must have sufficient horizontal 
expanse there is a natural map from  $\partial\bold H^2$ to itself.
The intuition behind the proof is that 
slicing $\h2r$ vertically gives Euclidean
planes, which have much less area than hyperbolic planes, 
so the image
of a horizontal hyperbolic plane under a quasi-isometry of $\h2r$ 
cannot be contained in a vertical slice.
The proof formalizes this idea by showing that if the image of vertical
geodesics, from a maximal family of geodesics all within a vertical cylinder 
all at least a certain distance apart, were to quasicross a vertical 
cross-section, then the number of quasicrossing points, all at least
a certain distance apart, exceeds the number possible in a Euclidean disk.

Key to Section 5 is the notion of semilocal growth of a finitely 
generated subgroup of a Lie group.  In particular, Theorem 5.19 
states that any finitely generated, non-elementary, non-discrete subgroup
of $PSL(2,\bold R)$ has superpolynomial semilocal growth. The semilocal
growth of a subgroup is strictly larger than its local growth, a notion
due to Carri\`ere [2]. In fact, all of the results of Section 5 can
be obtained replacing semilocal growth with local growth, but the notion
of semilocal growth seems more natural in this setting. It perhaps 
shows the importance of local growth and Theorem 5.19, which had
been previously obtained by Carri\`ere and Ghys [3], that they arose
independently in different contexts.

I would like to thank my advisor, Geoffrey Mess, for his guidance
throughout this project and for his enthusiasm for the subject. Thanks 
also go to Mladen Bestvina, Francis Bonahon and
Bob Edwards for many helpful conversations, and to 
Etienne Ghys for suggesting the proof of Lemma 5.8. I am grateful
to the superb anonymous referee whose careful reading and detailed
comments led to a significantly more elegant and succinct exposition.
Finally warmest thanks to my friends and family for their generous support, 
encouragement, and interest throughout this endeavor.

\head\newhead. Background and Definitions
\endhead

Two metric spaces are quasi-isometric if there is a relation
between them which, except locally, does not increase or 
decrease distances too much and has the property that
every point in each space is close to
a point which is related to a point in the other space.
\proclaim{Definition \newitem} A map $\psi:X\to Y$ is a 
$(\lambda,\epsilon,\delta)$--quasi-isometry
if it satisfies
\roster
\item (Lipschitz-in-the-large) $\frac{1}{\lambda}d(x_1,x_2)-\epsilon \le d(\psi(x_1),\psi(x_2))\le \lambda d(x_1,x_2)+\epsilon$, and
\item (Almost surjectivity) For all $y\in Y$, there 
exists $y^\prime\in Y$ with $d(y,y') \leq \delta$ such that, for some $x\in X$, $\psi(x)=y^\prime$.
\endroster
\endproclaim

There is some discrepancy in the literature as to whether quasi-isometries
are required to be almost surjective or not. Throughout this paper, when
we omit the $\delta$ in the definition of a quasi-isometry, we will mean
that we are not requiring the almost surjectivity condition, or in other 
words, the map is a quasi-isometry in the sense given above only onto the
range of the map. 

Two metric spaces are quasi-isometric if there exist quasi-isometries
between them that are almost inverses of each other. This idea is
formalized in the following definition.

\proclaim{Definition \newitem} Two metric spaces $X$ and $Y$ are 
$(\lambda,\epsilon,\delta)$--quasi-isometric if there exist 
$(\lambda,\epsilon,\delta)$--quasi-isometries $\psi:X\to Y$ 
and $\omega:Y\to X$ such that for some $\kappa$ 
\roster
\item $d(x,\omega(\psi(x)))\le \kappa$ for all $x\in X$, and
\item $d(y,\psi(\omega(y)))\le \kappa$ for all $y\in Y$.
\endroster
\endproclaim

Any group $\Gamma$ with generating set $S$, where $S$ contains the
inverse of any element in $S$, has a natural metric, the word metric $d_S$,
given by
	$$d_S(g,h)=\min\{n|ga_1a_2\dots a_n=h {\text { with }} a_i\in S\}.$$
Word metrics coming from two
different finite generating sets $S$ and $S'$  are equivalent in 
the sense that there is a constant $\lambda$ such that
	$$\frac{1}{\lambda}d_S(x,y)\leq d_{S'}(x,y)\leq \lambda d_S(x,y).$$
So up to quasi-isometry there is a natural metric on a finitely 
generated group.

\proclaim{Note \newitem}Whenever we say ``a group is 
quasi-isometric'' implicitly we will be talking about a finitely generated group
endowed with one of the equivalent word metrics.
\endproclaim

\head\newhead. A quasiaction of $\Gamma$ on $\h2r$
\endhead

Throughout this paper $\Gamma$ will denote a group 
$(\lambda,\epsilon,\delta)$--quasi-isometric
to $X$ via the word metric on $\Gamma$ coming from a finite
generating set. 

\proclaim{Definition \newitem} A 
$(\lambda_0,\epsilon_0,\delta_0,\kappa)$--quasiaction of a group $\Gamma$ on
a metric space $X$ consists of a family of 
$(\lambda_0,\epsilon_0,\delta_0)$--quasi-isometries 
$\phi_u:X\to X$ such that $d(\phi_u\phi_{u'}(x),\phi_{uu'}(x))\leq \kappa$ 
for all $u,u'\in\Gamma$ and all $x\in X$.
\endproclaim

Note that were $\kappa$ zero, we would have an action of $\Gamma$ on $X.$ Also
it's worth pointing out that the $\phi_u$'s are not necessarily
homeomorphisms. 

\proclaim{Observation\newitem} A group $\Gamma$ quasi-isometric to
a $\h2r$ has a natural quasiaction on $\h2r$, 
with $\kappa = \lambda\delta + \epsilon$,
given by $\phi_u = \psi \circ u \circ \omega$, where $\psi$ and 
$\omega$ are the 
$(\lambda, \epsilon, \delta)$--quasi-isometries
from the definition of $\h2r$ and $\Gamma$ being quasi-isometric, 
and $u\in \Gamma$ acts on $\Gamma$ by left multiplication.
\endproclaim
\name\psaction

\head\newhead. Constructing an action of $\Gamma$ on $\partial\bold H^2$. 
\endhead

Throughout this section ``vertical'' will refer to the $\bold R$
coordinates and ``horizontal'' to the $\bold H^2$ coordinates.
Let $\pi$ be the projection of $\h2r$ onto $\bold H^2\times\{0\}$.
In order to obtain an action of $\Gamma$ on $\partial\bold H^2$ of $\bold H^2$
we want to show that each $\pi\circ\phi_u$ restricted to the horizontal plane
$\bold H^2\times\{0\}$ is a quasi-isometry. The idea is to show 
that  each $\phi_u$ must preserve horizontal and vertical
 in some rough sense.
The proof exploits the fact that there is a lot more room in a disk in
the hyperbolic plane than in a disk in the Euclidean plane. We begin with
two estimates exhibiting this difference.

\proclaim{Estimate \newitem} An upper bound for the number of disks of radius $r$
in a disk of radius $R$ all at least $2s$ apart in Euclidean
space is given by

          $$\left(\frac{(R + s)}{(r + s)}\right)^2.$$ 

\endproclaim   
\name\Rest

Recall that the area of a hyperbolic disk of radius $R$ is given by
$2\pi(\cosh R-1).$

\proclaim{Estimate \newitem}A lower bound for a maximal number of disks 
of radius $r$ in a disk of radius $R$ all at least $2s$ apart in hyperbolic
space is given by

	$$\frac{e^R-2}{2(\cosh(2(r+s))-1)}.$$ 

\endproclaim
\name\Hest

By ``the horizontal distance between points $x$ and $y$,'' I will mean 
the distance in $\bold H^2\times\{0\}$ between $\pi(x)$
and $\pi(y)$, and by ``vertical distance,'' the distance between the 
$\bold R$--coordinates.

\proclaim{Proposition \newitem}Let $\phi:\h2r\to\h2r$ be a 
$(\lambda,\epsilon)$--quasi-isometry. Let $C_{0}$ be a vertical solid cylinder
of radius $r$, and let $p_{0}$ and $q_{0}$ be points in $C_{0}$ such that
$\phi(p_{0})$ and $\phi(q_{0})$ have vertical distance no greater than $h_0$.
Then $p_{0}$ and $q_{0}$ must be no farther apart than 
	$$ch_0+c$$
where $c$ is a constant depending only on $\lambda$, $\epsilon$ and $r$.
\endproclaim
\name\mainprop

\demo{Proof}Let $d(p_{0},q_{0})=L$. Let $w$ denote the horizontal distance between 
$\phi(p_{0})$ and $\phi(q_{0})$ and $h$ the vertical. Now

	$$\align
	L/\lambda-\epsilon   &\leq d(\phi(p_{0}),\phi(q_{0})) \\
		&=\sqrt{w^2+h^2} \\
		&\leq w+h_0.	
	\endalign$$ 
So

	$$w \geq L/\lambda -\epsilon -h_0.$$
Project $\phi(p_{0})$  and $\phi(q_{0})$ onto $\bold H^2\times\{0\}$.
Let $\alpha$ be the geodesic which intersects the
geodesic between the projected images of $\phi(p_{0})$ and $\phi(q_{0})$ 
perpendicularly at the point halfway between them. 

For every vertical solid cylinder $C$ of radius $r$ whose central
axis is distance $D$ from the central axis of $C_{0}$, there are points
$p$ and $q$ that are the translation of $p_{0}$ and $q_{0}$ under 
the translation that takes the central axis of $C_{0}$ to the central axis of
$C$ and preserves the vertical height of each point.  
Let $\gamma$ be the geodesic between $p$ and $q$. 
The points $\phi(p)$ and $\phi(q)$ are on opposite sides of 
$\{\alpha\}\times \bold R$, if
	$$d(\phi(p_{0}),\phi(p))\leq 1/2w$$  
and 
	$$d(\phi(q_{0}),\phi(q))\leq 1/2w.$$

Let $R$ be the maximal distance $D$ for which these inequalities hold:
	$$R=\frac{1}{2\lambda}(L/\lambda-3\epsilon-h_0).$$ 
Let $\Cal C$ be a maximal set of cylinders of radius $r$ all $s$ apart within
the cylinder of radius $R$ about $C_{0}$, where $s = \lambda(1 + 2\epsilon)$.
For each cylinder $C$ in  $\Cal C$, let $x$ be a point in 
$\{\alpha\}\times \bold R$ that is as close as possible to $\phi(\gamma)$.
Such points are called quasicrossing-points of  $\phi(\gamma)$ with
$\{\alpha\}\times \bold R$ since $d(x,\phi(\gamma)) \leq \epsilon/2$. Two
quasicrossing-points $x$ and $x'$ associated with cylinders $C$ and
$C'$ in  $\Cal C$ must be at least $1$ apart as
$$\align
     d(x,x') 
          &\geq d(y,y') - d(x,y) - d(x',y')\\
          &\geq s/\lambda - \epsilon - \epsilon/2  - \epsilon/2 =1,
\endalign$$
where $y$ and $y'$ are points 
on $\phi(\gamma)$ and $\phi(\gamma')$ 
within $\epsilon/2$ of $x$ and $x'$ respectively.

To see that
$$d(x,x') \leq \lambda(L + R) + 2\epsilon,$$
let $z$ and $z'$ be
points on $\gamma$ such that $\phi(z) = y$ and $\phi(z')=y'$. 
Since all points on  $\gamma$ and  $\gamma'$ lie within a cylinder
of radius $R$ and height $L$, $d(z,z') \leq L + R$. So
$$\align
d(x,x') &\leq d(y,y') + d(x,y) + d(x',y')\\
        &\leq \lambda d(z,z') + \epsilon + \epsilon\\
        &\leq \lambda(L +R) + 2\epsilon.
\endalign$$

By Estimate \Hest\ there are at least 
$\frac{e^R-2}{2(\cosh(2(r+s/2))-1)}$ 
vertical cylinders of radius $r$ in $\Cal C$.
Let $X$ be a set of quasicrossing-points, one for each cylinder
$C$ in $\Cal C$. 
By estimate 1 there are at most 
          $$(2\lambda(L + R) +4\epsilon +1)^2$$
such points. Thus 
$$\frac{e^R-2}{2(\cosh(2(r+s/2))-1)} \leq (2\lambda(L + R) +4\epsilon +1)^2.$$
Solving for $L$ in the definition of $R$ gives
	$$L=2\lambda^2R+3\lambda \epsilon +\lambda h_0.$$
Substituting this expression for $L$ in the above inequality shows that
$e^R$ must be at most some quadratic function of $\lambda$, $\epsilon$,
$r$ and $h_0$. Furthermore 
$$e^R \leq b_2 R^2 + b_1 R + b_0,$$
where the coefficients $b_0$, $b_1$ and $b_2$ are polynomial
functions of $h_0$ of degree no more than 2. 
The first four terms of the Taylor expansion of $e^R$ tell us that
$R$ must be no more than $b_2 + 2b_1  + 6b_0$. So $R$ is no more than
some quadratic function of $h_0$. The definition of $R$ together with
this bound give a similar bound on $L$,
$$L \leq c_2 h_0^2 + c_1 h_0 +c_0,$$
where the coefficients depend only on  $\lambda$, $\epsilon$, and $r$.
We can get a linear bound on $L$ with respect to $h_0$ as follows.
Let $c = c_2 + c_1 + c_0$. A geodesic between $\phi(p_{0})$ 
and $\phi(q_{0})$ of vertical distance $h\leq h_0$ can be split up into
no more than $h_0+1$ pieces of vertical distance no greater than 1. 
Using the above
bound on each piece and adding them together, we know that $p_{0}$ and
$q_{0}$ can be no farther apart than $c(h_0+1)$.
\qed\enddemo

\proclaim{Remark \newitem}For future reference note that $c=aS$ where 
$S = 2(\cosh(2r + \lambda(1 + 2\epsilon)) -1)$
and  $a$ is some constant depending only on $\lambda$ and $\epsilon$. 
\endproclaim
\name\lininS

\proclaim{Corollary \newitem}Let $\phi_u:\h2r\to\h2r$ be a
($\lambda,\epsilon$)--quasi-isometry coming from the quasiaction of
$\Gamma$ on $\h2r$.
Let $x$ and $y$ be points
$D$ apart in some horizontal $\bold H^2$. Then the 
horizontal distance $l$ between $\phi_u(x)$ and $\phi_u(y)$ must be
at least
	$$l\geq\ln((\frac{\frac{D}{\lambda}-\epsilon}{(a(2\kappa + 2\delta+1)}-1)-\lambda (2\epsilon+1) + 1$$ 
where $\kappa$ is the quasiaction constant of Observation \psaction, and
$a$ depends only on $\lambda$ and $\epsilon$.
\endproclaim
\name\horizontal
\demo{Proof}Since $d(\phi_{u^{-1}}\phi_u(x),x)\leq d(\phi_{u^{-1}}\phi_u(x), \phi_e(x)) + d(\phi_e(x), x) \leq \kappa + \delta$ and
similarly for $y$, the vertical distance between
$\phi_{u^{-1}}(\phi_u(x))$ and $\phi_{u^{-1}}(\phi_u(y))$ 
must be less than $2(\kappa + \delta)$.
Hence by Proposition \mainprop, $d(\phi_u(x),\phi_u(y))\leq 2c(\kappa +\delta)+c$. By Remark \lininS,
	$$d(\phi_u(x),\phi_u(y))\leq 2a\cosh(2r+\lambda (2\epsilon+1))-1)(2(\kappa+\delta) + 1),$$
where $r$ is the radius
of a vertical cylinder containing $\phi_u(x)$ and $\phi_u(y)$, so can be taken
to be half the horizontal distance $l$ between $\phi_u(x)$ and $\phi_u(y)$. 
It is $l$ we are trying to bound from below. So
	$$\frac{D}{\lambda}-\epsilon\leq d(\phi_u(x),\phi_u(y))\leq (2a\cosh(l+\lambda  (2\epsilon+1) -1)(2(\kappa + \delta) + 1),$$
from which the result follows.
\qed\enddemo

\proclaim{Proposition \newitem}$\pi\circ\phi:\bold H^2\times\{0\}\to
\bold H^2\times\{0\}$ is
a ($\lambda',\epsilon'$)--quasi-isometry where $\lambda'=\max\{\lambda,D'\}$ and 
$\epsilon'=\max\{\epsilon,1\}$ taking $D'=\lambda a(2\kappa + 2\delta +1)(e^{\lambda (2\epsilon+1)}+1) + \epsilon).$                 
\endproclaim
\name\psiso

\demo{Proof} The projection $\pi$ is distance 
nonincreasing, so for any two points
$x$ and $y$ in $\bold H^2\times\{0\}$ we know 
	
	$$d(\pi\circ\phi(x),\pi\circ\phi(y))
			\leq \lambda d(x,y)+\epsilon.$$
We need to show that $\pi\circ\phi$ does not decrease distances too
much. By Corollary\horizontal\, the images under $\pi \circ \phi$
of any two points at least distance
$D'$ apart cannot be closer than 1 unit apart. Let 
$L=d(\pi \circ \phi(x),\pi \circ \phi(y))$, and 
let $\gamma$ be the geodesic between $\pi \circ \phi(x)$ and 
$\pi \circ \phi(y)$.
Let the $x_i$ be $N\leq L+1$ points along $\gamma$ such that
successive $x_i$ are within distance 1 of each other and 
$x_0 = \pi \circ \phi(x)$ and $x_N = \pi \circ \phi(y)$. So 
successive quasi-preimages of the $x_i$ must be no more
than $D'$ apart. Hence
	$$d(x,y)\leq D'(L+1)=D'd(\pi \circ \phi(x),\pi \circ \phi(y))+D'.$$
Therefore,
	$$d(\pi \circ \phi(x),\pi \circ \phi(y))\geq\frac{1}{D'}d(x,y)-1.$$
\qed\enddemo
 
A quasisymmetric map of $S^1$ viewed as the boundary
of the Poincar\'e model for the hyperbolic plane is a map that extends to a 
quasiconformal map of $\bold H^2$.
If we view $\bold H^2$ as the upper half plane of $\bold C$,
Beurling and Ahlfors [1] showed that $f$ is quasisymmetric
and fixes the point at infinity exactly when there exists
a constant $\epsilon$ such that
	$$1/\epsilon\leq \frac{f(x+t)-f(x)}{f(x)-f(x-t)}\leq \epsilon.$$
It is well-known (see for example [9][16][26]) that any quasi-isometry
of $\bold H^2$ extends to a map on $\partial\bold H^2$ that is
quasisymmetric and the quasisymmetry constant only depends
on the quasi-isometry constants. The central idea is that the image
of any geodesic under a quasi-isometry lies in the $B$
neighborhood of a geodesic, where the constant $B$ depends only on
$\lambda$ and $\epsilon$. Although we only had a quasiaction of $\Gamma$
on $\bold H^2\times\bold R$ and also on $\bold H^2\times \{0\}$,
we get a true action on $\partial\bold H^2$ for the following reasons.
Since the images of a geodesic in $\bold H^2\times\{0\}$ under
$\pi\circ\phi_{uu'}$ and $(\pi\circ\phi_u)\circ(\pi\circ\phi_{u'})$
are within bounded distance of each other, they must be within bounded
distance of the same geodesic. To see where a point $x$ in $\partial\bold H^2$
goes under a quasi-isometry $\phi$, take a geodesic ray 
$\gamma(t)$ with $x$ as its
endpoint. The subset $\phi(\gamma)$ is within a bounded distance of some
geodesic ray.  Map $x$ to the appropriate endpoint of this geodesic ray.
Since this geodesic ray is the same for both maps, both maps must give the
same action on the boundary.

\proclaim{Corollary \newitem} There is a canonical surjective homomorphism
$\Xi:\Gamma\to F$ where $F$ is the uniformly quasisymmetric
group consisting of the boundary values of the $\pi\circ\phi_u$'s.
\endproclaim

\head\newhead. The action of $\Gamma$ on $\partial\bold H^2$ is 
conjugate by a quasisymmetric map to an action by a M\" obius group.
\endhead

Let $F^+$ be the group of orientation  preserving elements of $F$.
This section shows that $F^+$ may be conjugated by a quasisymmetric
map to a M\"obius group $G$.
The results of this section were previously known. For example,
Lemma 4.3 is a special case of Theorem 9 in [25], and
Lemma 4.4 follows trivially from results in Section
C of [35].

For $F^+$ not discrete in the topology of pointwise convergence
(a case we will rule out in section 5), we may use the following
theorem of Hinkkanen.

\proclaim{Theorem \newitem (Hinkkanen) [21]}Let $G$ be a 
uniformly quasisymmetric group containing a sequence of distinct 
elements that tend to the identity pointwise. Then $G$ is
a quasisymmetric conjugate of a M\"obius group.
\endproclaim

Throughout the rest of the section we will assume $F^+$ is
discrete. Any discrete uniformly quasisymmetric group of
orientation preserving homeomorphisms is a convergence group,
so the following theorem holds for $F^+$.

\proclaim{Theorem \newitem (Casson--Jungreis [6], Gabai [14], Tukia [36])}
$G$ is a convergence group if and only if $G$ is conjugate
in Homeo($S^1$) to the restriction of a Fuchsian group.
\endproclaim

This theorem was proved by Gabai [14], and independently by
Casson--Jungreis [6], building on work of Tukia [36].
The rest of this section is devoted to showing that the map
conjugating $F^+$ to a M\"obius group may be taken to be
quasisymmetric.

\proclaim{Lemma \newitem} Say $\Theta:\Gamma\to\bold H^2$ and 
$\Psi:\Gamma\to\bold H^2$ are two maps which induce a 
quasi-isometry between $\Gamma$ and $\bold H^2$. Then the 
quasisymmetric maps they induce on $\partial\bold H^2$ are conjugate 
by a quasisymmetric map.
\endproclaim
\name\conjBvalues
\demo{Proof} Define $H:\Theta(\Gamma)\to\Psi(\Gamma)$ to be
the composition of $\Psi$ with a quasi-inverse of $\Theta$. 
As we saw in the previous section 
$H$ can be extended to a quasi-isometry of $\bold H^2$ 
to itself, which induces a quasisymmetric map $h$ on 
$\partial\bold H^2$. By examining the construction, we see that $h$
conjugates the map induced on $\partial\bold H^2$ by $\Psi$ to the
one induced by $\Theta$.
\qed\enddemo

\proclaim{Lemma \newitem} For any $(\lambda,\epsilon)$--quasi-isometry 
$f:\bold H^2\to\bold H^2$ that has the same boundary values as some 
isometry $g:\bold H^2\to\bold H^2$, there is a constant $L$ such that 
$d(f(x),g(x))\leq L$ for all $x\in \bold H^2$, where $L$ is
dependent only on $\lambda $ and $\epsilon$.
\endproclaim
\name\psclosetoiso

\demo{Proof}Let $y_0,y_1,y_2$ be vertices of an ideal triangle such 
that $x$ is within the
same distance, say $P$, of each of the geodesics connecting the
three points. Denote by $Y_i$ the geodesic connecting $y_i$ and
$y_{(i+1)mod 3}$.
As was explained in the paragraph following the proof of 
Lemma \psiso\, , $f(Y_i)$ is a curve that remains within
distance $D$ (depending only on $\lambda$ and $\epsilon$) 
of $g(Y_i)$. The point
$f(x)$ is within distance $\lambda P+\epsilon$ of each $f(Y_i)$, 
so must be within $\lambda P+\epsilon+D$ of each $g(Y_i)$. 
The intersection of these regions has bounded
diameter, say $L$. Since $g(x)$ is certainly within this region,
$f(x)$ must be within $L$ of $g(x)$.
\qed\enddemo

\proclaim{Lemma \newitem}Any finitely generated, discrete, uniformly
$(\lambda ,\epsilon)$-quasisymmetric group $G$ acting on $S^1$, that induces
a cocompact action on the space of triples $T$, is 
quasi-isometric to $\bold H^2$.
\endproclaim
\name\ccptontriples

\demo{Proof}To any ordered triple of distinct points in $S^1$ we associate
a point in hyperbolic space as follows. Connect the first two points 
by a geodesic and drop a perpendicular from the third. The intersection
will be the point associated to the triple. Choose  a triple $t_0$ 
whose associated point is $x_0$.
Let $f:G\to\bold H^2$ be the map sending an element
$a\in G$ to the point associated to the triple $a(t_0)$. We will
show that $f$ is a quasi--isometry. For any element $a\in G$, choose a 
$(\lambda,\epsilon)$--quasi-isometry $\eta_a:\bold H^2\to\bold H^2$ with boundary 
values $a$. Note $\eta_a(f(a'))$ must be within $2B$ of $f(aa')$ (where $B$
was defined  in the paragraph following Proposition \psiso), since 
$\eta_a(f(a'))$ must be within $B$ of each of the geodesics used to construct
$f(aa')$. 

Since $G$ is cocompact on the space of triples, there is a constant
$E$ such that every point in $\bold H^2$ is within $E$ of some point in
the  image of $f$. Let $J'$ be a finite generating set for $G$.
Enlarge $J'$  to $J$ by including all elements $a$ such that
$d(\eta_a(x_0),x_0) \leq 3\lambda E+\epsilon + 6B$. 
By the discreteness of $G$, there 
are only finitely many
such elements.  Let $M=\max_{b\in J} d(f(b),x_0)$. We are now ready to
check that $f$ is a quasi-isometry. If $d(a,a')=m$, then for
some $b_i\in J$, $a=a'b_1b_2\dots b_m$. So,
	$$\align
	d(f(a),f(a'))&\leq d(f(a'b_1\dots b_m),f(a'b_1\dots b_{m-1}))\\
		     &\qquad+d(f(a'b_1\dots b_{m-1}),f(a'b_1\dots b_{m-2}))
			+\dots +d(f(a'b_1),f(a')) \\
		&\leq d(\eta_{a'b_1\dots b_{m-1}}f(b_m),\eta_{a'b_1\dots b_{m-1}}
f(e))+\dots +d(\eta_{a'}f(b_1),\eta_{a'}f(e)) + 4mB \\
		&\leq \lambda d(f(b_m),f(e))+\epsilon+\dots +\lambda d(f(b_1),f(e))+\epsilon + 4mB\\
		&\leq(\lambda M+\epsilon +4B)m \\
		&=(\lambda M+\epsilon+4B)d(a,a').
	\endalign$$

To get the lower bound, we compute as follows. Let $d(f(a),f(a'))=L$.
Divide the geodesic between $f(a)$ and $f(a')$ into $[L/E]+1$ segments
of length no more than $E$, with endpoints $x_0,x_1,\dots,x_N$, where
$f(a)=x_0$ and $f(a')=x_N$. Each $x_i$ is within $E$ of some point $f(a_i)$. 
By construction, for every $i$, $d(f(a_{i-1}),f(a_{i}))\leq 3E$. We wish
to show that there is a generator taking $a_i$ to $a_{i + 1}$ for all $i$.
Recall that $d(\eta_a(f(a'), f(aa')) \leq 2B$ and $x_0 = f(e)$. Note that
$$\align
d(x_0, f(a_{i-1}^{-1}a_i)) &\leq d(x_0, \eta_{a_{i-1}^{-1}}(f(a_{i-1})))
         + d(\eta_{a_{i-1}^{-1}}(f(a_{i-1})),\eta_{a_{i-1}^{-1}}(f(a_{i})))\\
            &\qquad + d(\eta_{a_{i-1}^{-1}}(f(a_{i})), f(a_{i-1}^{-1}a_i))\\
            &\leq 2B + 3\lambda E +\epsilon +2B.
\endalign$$
Furthermore,
$$\align
d(\eta_{a_{i-1}^{-1}a_i}(x_0), x_0) &\leq d(x_0, f(a_{i-1}^{-1}a_i))) 
          + d(f(a_{i-1}^{-1}a_i), \eta_{a_{i-1}^{-1}a_i}(x_0)\\
          &\leq 3\lambda E + \epsilon + 6B.
\endalign$$
So for all $i$, $a_{i-1}^{-1}a_i = b_i$ is in the generating set for $G$.

Since
	$$a'=a_N=a_{N-1}b_N=\dots=ab_1b_2\dots b_N,$$
we have
	$$\align
	d(a,a')&\leq N\\
		&= [L/E]+1\\
		&\leq \frac1E d(f(a),f(a'))+1.
	\endalign$$
\qed\enddemo

\proclaim{Observation \newitem}Section 2 together with 
the last paragraph of section 3 shows how
a quasi-isometry of a group $G$ with the hyperbolic
plane induces an action on $\partial\bold H^2$. It is easy to check
that the action of $F^+$ on $\partial\bold H^2$ that we get from the 
quasi-isometry of $F^+$ with $\bold H^2$ given by the Lemma \ccptontriples\,
is the same action of $F^+$ on $\partial\bold H^2$ that we started with.
\endproclaim
\name\sameaction

\proclaim{Proposition \newitem}$F^+$ is conjugate to a M\" obius group by a
quasisymmetric homeomorphism $h:S^1\to S^1$.
\endproclaim

\demo{Proof}For $F^+$ not discrete this is a  result 
of Hinkkanen [21].

Gabai [G], or Casson--Jungreis [6], together
with Tukia [36] show that any 
discrete convergence group (in particular, any quasisymmetric group)
is conjugate to a finite extension of a Fuchsian group 
by some homeomorphism $h_1$ of $S^1$. 
Our $F^+$ acts cocompactly on the space of triples, so $G^+=h_1F^+h_1^{-1}$
must also. Thus $G^+$ is a Fuchsian group acting  discretely and 
cocompactly on the space of triples, so when we extend the action 
to $\bold H^2$ we get a discrete cocompact group of hyperbolic isometries.
This gives us a quasi-isometry of $G^+$ with $\bold H^2$.
We get a different quasi-isometry $\Psi:G\to\Psi(G)$ by using
the isometry of $G^+$ with $F^+$ and applying Lemma \ccptontriples\
to $F^+$.
Applying Lemma \conjBvalues\ to these two quasi-isometries,
which have the same boundary values by Observation \sameaction\,,
yields the desired result.

\qed\enddemo

Let $\Phi:\Gamma\to G$ be the homomorphism we have obtained where
$G$ is either the M\"obius group $G^+$ or a $\bold Z/{2\bold Z}$ extension
of $G^+$, depending on whether $F$ contained orientation reversing 
elements or not.

\head\newhead. The discreteness of the image of $\Phi$
\endhead

Our goal now is to show that $G$ is discrete.  The idea is that 
any non-elementary M\"obius group which is not discrete has 
many more small elements than
$G$. More precisely, choose a left-invariant metric on $PSL(2,\bold R)$.
For any $\epsilon >0$, denote by $N_\epsilon$ 
the $\epsilon$-neighborhood of the identity in  
$PSL(2,\bold R)$. For any finitely generated subgroup 
$H$ of $PSL(2,\bold R)$ , let $H^n_\epsilon$ denote
the set of $h_u \in N_\epsilon$ such that $u$ has word length 
less than or equal to $n$. Our claim is that for some $\epsilon$
(specified later), $|G^n_\epsilon|$,
the number of elements in $G^n_\epsilon$, grows linearly 
with $n$, while for any finitely generated nonelementary group $H$ 
that is not discrete $|H^n_\epsilon|$ grows exponentially in $n$.

\proclaim{Definition \newitem} Let $G$ be a finitely generated
group with word metric $d$ imbedded in another group $L$ with metric
$\rho$. Then the semilocal growth of $G$ in $L$ is defined to be the growth
rate of the number of elements in 
	$$G^n_\epsilon=\{g\in G| d(g,e)\leq n, \rho(g,e)\leq\epsilon\}.$$
\endproclaim

\proclaim{Note \newitem} This notion of local growth is strictly
bigger than Carri\`ere's notion of local growth [2], which counts
only elements $g$ such that the subwords of increasing length making up $g$
are all in the $\epsilon$-neighborhood of the identity. More, precisely
it counts only $g=a_{i_1}a_{i_2}\dots a_{i_m}$ where 
$g_j = a_{i_1}a_{i_2}\dots a_{i_j}\in N_\epsilon$ for all $j \leq m$, where
the $a_k$ are generators for $G$. All of the results hold equally well
for local growth as for semilocal growth, but as local growth 
imposed an extra, unnecessary condition on the elements we are counting, we 
will prove the results for semilocal growth.
\endproclaim

\subheading {$G$ has linear semilocal growth in $PSL(2,\bold R)$}

Recall that $\psi:\Gamma\to \h2r$ is the quasi-isometry that came from the
definition of $\Gamma$ and $\h2r$ being quasi-isometric.

\proclaim{Observation \newitem} The number of $u\in \Gamma$ such that
the images of $z_0=\psi(e)$ under $\phi_u$
lies in a vertical cylinder $K$ of height 1 and radius $R$
centered about the vertical geodesic through $z_0$ is bounded by 
some constant $N$ depending only on $\lambda$, $\epsilon$ and $R$.
\endproclaim
\name\noofimages

\proclaim{Observation \newitem}The number of $u\in \Gamma$ 
of word length $n$ such that
$\phi_u(z_0)$ lies in $C_R$, the vertical cylinder of radius
$R$ centered about $z_0$, is bounded by $An+B$ where
$A=2\lambda N$ and $B=(2\epsilon+1)N$.
\endproclaim
\name\noofwords

\proclaim{Theorem \newitem} For any $\epsilon$, $|G_\epsilon^n|$ grows 
linearly with $n$. More explicitly, there exist constants 
$A$ and $B$ such that the number of elements in $G_\epsilon^n$ 
is less than or equal to $An+B$.
\endproclaim
\name\lingrowth

\demo{Proof} For any $u\in\Gamma$ let
$g_u = \Phi(u)$. 
There is some constant $r_\epsilon$ such that 
any $g_u$ within $\epsilon$ of the identity moves $x_0$ no
more than $r_\epsilon$. Let $f_u$ be the map whose conjugate under $h$
is $g_u$. By the Lemma \psclosetoiso\,, $f_u$ must not 
move $x_0$ by more than $L+r_\epsilon$. The map
$f_u$ came from the projection of the image of the horizontal
$\bold H^2$ containing $\psi(e)$ under $\phi_u$. Since $\psi(e)$
projects to $x_0$, we conclude that $\phi_u$ must send $\psi(e)$
to some point in the vertical cylinder of radius $L+r_\epsilon$
centered about $\psi(e)$. Taking $R= L + r_\epsilon$ in 
Observation \noofwords\,, 
the number of such $u$ of word length less than or equal 
to $n$ is bounded by $An+B$.
\qed\enddemo

\proclaim{Lemma \newitem}The kernel of $\Phi$ is either finite
or contains an element of infinite order.
\endproclaim

\demo{Proof}Passing if necessary to an index 2 subgroup, we 
may assume that all elements $u$
of $\ker \Phi$ are end--preserving in the sense that the 
image under $\phi_u$ of
a sequence of points whose $\bold R$ components tend to $+\infty$
(resp. $-\infty$) also have $\bold R$ components which tend to 
$+\infty$ (resp. $-\infty$). 
Our goal will be to show that if an element $u\in\ker\Phi$
has finite order, then $\psi(u)$ and $\psi(e)$ are within a vertical
distance of $\lambda c+\epsilon$ of each other, where $c$ is the constant of
Proposition \mainprop. This would imply that there are only finitely
many elements of finite order, since the image under $\Psi$ of any element
of finite order in the kernel would have to be in the 
vertical cylinder $C_{0}$ of radius $L$ centered at $z_0$
bounded above and below by 
$\lambda (c\epsilon/2 + c)  \epsilon + 2(\lambda\delta + \epsilon)$.
But by Observation \noofimages,
the number of such points is finite.

Say $u$ has finite order $n$ and the vertical distance between
$\psi(u)$ and $\psi(e)$ is greater than $\lambda c+\epsilon$. We assume without
loss of generality that $\psi(u)$ is above $\psi(e)$.
Then for some $k$,
$\psi(u^k)$ is above $\psi(u^{k+1})$. Draw a geodesic segment
between these two points and then continue it vertically upward 
from $\psi(u^k)$ and vertically downward from $\psi(u^{k+1})$.
Call this curve $\gamma$. Under $\phi_{u^{n-k}}$, $\psi(u^k)$
is sent to $\psi(e)$ and $\psi(u^{k+1})$ is sent to $\psi(u)$. Since
we are assuming $\phi_u$ is end--preserving, it follows that there 
is some point $p$ on $\gamma$ above $\psi(u^k)$, such that 
$\phi_u(p)$ is within $\epsilon/2$ of being at the same height 
as $\psi(u)$. From the previous
section we know that this means that $p$ and $\psi(u^{k+1})$ must
be no farther than $c\epsilon/2 + c$ apart vertically. Thus,
$$d(\psi(u^{k+1}),\psi(u^k))\leq d(\psi(u^{k+1}),p)\leq c\epsilon/2 + c.$$
Therefore,
	$$\align
	d(\psi(e),\psi(u))&\leq d(\phi_{u^{n-k}}\psi(u^{k+1}),\phi_
			{u^{n-k}}\psi(u^k)) + 2(\lambda\delta + \epsilon)\\
		&\leq \lambda d(\psi(u^{k+1}),\psi(u^k))+\epsilon  + 2(\lambda\delta + \epsilon)\\
		&\leq \lambda (c\epsilon/2 + c)  \epsilon + 2(\lambda\delta + \epsilon).
	\endalign$$
\qed\enddemo
\proclaim{Theorem \newitem}The kernel is either finite or 
quasi-isometric to $\bold R$.
\endproclaim

\demo{Proof}Suppose the kernel is infinite. Let $u\in\ker\Phi$ 
be an element of infinite order, whose existence is guaranteed 
by the previous Lemma. We may assume without loss
of generality that $\phi_u$ is end--preserving and that $u$ was in our 
generating set for $\Gamma$ to begin with. Since $u$ is of infinite
order, there are infinitely many $\psi(u^i)$. Since all $\psi(u^i)$
are in $C_{0}$, the vertical distance between
$\psi(u^n)$ and $z_0$ must be greater than $c$, for $n$ 
sufficiently large. Let $u^n=v$. Without loss
of generality we may assume that $\psi(v)$ is above $z_0$. By an
argument similar to that of the previous Lemma, $\psi(v^{k+1})$ is
above $\psi(v^k)$. Let $d(\psi(e), \psi(v))=h$. Then for all $k$,
$d(\psi(v^{k+1}),\psi(v^k))\leq \lambda h+\epsilon$. Thus every point in $C_{0}$
is within $\frac12(\lambda h+\epsilon)+L$ of one of the $\psi(v^k)$.

  Say 
$w\in\ker\Phi$. Then $\psi(w)$ is within $\frac{1}{2}(\lambda h+\epsilon)+L$ of 
some $\psi(v^k)$. Thus,
	$$\align
      d(z_0,\psi(v^{-k}w))&\leq \lambda d(\phi_{v^k}(z_0),\phi_
{v^k}(\psi(v^{-k}w))+\epsilon\\
		&=\lambda d(\psi(v^k),\psi(w))+\epsilon\\
		&\leq \lambda (\frac{1}{2}(\lambda+\epsilon)+L)+\epsilon.
	\endalign$$
There are only finitely many such points so $\langle v\rangle$ is of
finite index in $\ker\Phi$. Thus $\ker\Phi$ is quasi-isometric
to $\bold Z$ which in turn is quasi-isometric to $\bold R$.
\qed\enddemo
\name\finorcqi
\proclaim{Lemma \newitem}The group $G$ is a non-elementary
subgroup of $PSL(2,\bold R)$.
\endproclaim

\demo{Proof}By construction $G$ does not fix a point or
preserve an axis in $\bold H^2$. Say $G$ fixes a point on 
$\partial\bold H^2$. Then $G$ would be solvable. Since 
$\ker\Phi$ is either finite or a finite extension of $\bold Z$,
$G$ solvable would imply that $\Gamma$ was amenable. But
$\Gamma$ cannot be amenable since it is quasi-isometric 
to $\bold H^2\times\bold R$.
\qed\enddemo
\name\nonelt

\proclaim{Theorem \newitem}The group $G$ is a non-elementary subgroup
of $PSL(2,\bold R)$ with the property that $|G^n_\epsilon|$ grows
linearly.
\endproclaim

\subheading{Non-discrete, finitely generated, non-elementary subgoups of 
$PSL(2,\bold R)$ have exponential semilocal growth}

In order to show that $|H_{\epsilon}^n|$ grows quickly we need
a tool for constructing small elements and another to make sure we
can construct enough. The Zassenhaus Lemma will perform the
first task while Tits's Theorem will do the second. We need the 
following Lemmas before we can apply the results. 

\proclaim{Lemma \newitem}$G$ is either discrete or its closure
is all of $PSL(2,\bold R)$.
\endproclaim
\name\fullclosure
\demo{Proof}By construction $G$ does not fix any point in $\bold H^2$ or
its boundary, and also does not preserve any axis. It is a well-known
fact (see for example [4], Theorem 4.4.7) that the only closed subgroups
of $PSL(2,\bold R)$ either fix one of the above, are discrete, or
are all of $PSL(2,\bold R)$.
\qed\enddemo

Let $H$ be a subgroup of $PSL(2,\bold R)$ that is not discrete and whose
closure is all of $PSL(2,\bold R)$.

\proclaim{Lemma \newitem}$H$ is not virtually solvable.
\endproclaim

\demo{Proof}The Lemma follows from the fact that a virtually solvable subgroup
has virtually solvable closure.
\qed\enddemo

\proclaim{Lemma \newitem}(Zassenhaus [37]) There is a constant $\epsilon_0$
such that for any $r\leq\epsilon_0$, if $f$ and $g$ are within distance       
$r$ from the identity, then $[f,g]$ is also.
\endproclaim

A proof of this Lemma may be found in Raghunathan's book [28].

We will use this Lemma to construct more small elements from a few
small elements. But we need a way to check that we  are indeed
constructing new and  different elements by conjugating.
Tits's alternative says that any subgroup of a linear group is either
virtually solvable or contains a free group on two generators. We 
need a little more; we need a free group generated by two small
generators in the Zassenhaus sense.

\proclaim{Note \newitem}Unknown to me at the time, Carri\`ere and
Ghys [3] had previously proved the existence 
of such a free group along similar lines.
\endproclaim

We will construct these elements using the following fact. 

\proclaim{Proposition \newitem}Any finitely gen\-er\-at\-ed 
non\-el\-e\-men\-tary  
sub\-group $H$ of the group $PSL(2,\bold R)$ that is not 
discrete must contain an elliptic element of infinite order.
\endproclaim

\demo{Proof}
By Selberg's Lemma, $H$ contains a finite index normal subgroup $N$
containing no non-trivial elements of finite order, from which it follows
that $H$ has no finite order elements of order greater than the index of
$N$ in $H$.
But $H$ is dense in $PSL(2,\bold R)$,
so there must be elements of $H$ which get close to finite order elliptics
of higher orders. The only possibilities are elliptics of infinite order,
so $H$ must contain at least one.
\qed\enddemo
 
Let $\alpha$ be an infinite order elliptic element lying 
in $H$ whose existence is guaranteed by
the previous Proposition. The eigenvalues of $\alpha$ all have
norm 1, but are not roots of unity, so they have infinite order in 
$k^*$, the subfield of $\bold C$ generated by the matrix entries and
eigenvalues of the generators for $H$. Thus we may use 
the following Lemma of Tits [32, 4.1], taking $t$ to be one of 
$\alpha$'s eigenvalues. 

\proclaim{Lemma \newitem}Let $k$ be a finitely generated field and let $t\in k$
be an element of infinite order. Then there exists a locally compact
field $k'$ endowed with an absolute value $\omega$ and a homomorphism
$\sigma:k\to k'$ such that $\omega(\sigma(t))\ne 1$. 
\endproclaim

Tits finds subsets $S$ of $H$ and 
$\Cal U$ of $SL(2,\bold R)\times SL(2,\bold R)$ with the property
that for any $s_1,s_2\in S$ with $s_1\ne s_2$ and $(s_1,s_2)\in\Cal U$
there exists a positive power $m$ such that $s_1^m$ and $s_2^m$
generate a free group. He constructs these sets as follows.

The set $S$ consists of elements $s\in H$ with eigenvalues
$\lambda_1$, $\lambda_2$ such that $\omega(\sigma(\lambda_i))\ne 1$ for
$i=1,2$. Let $\Cal U$ be the set of 
$(x,y)\in SL(2,\bold R)\times SL(2,\bold R)$ 
 such that $x$ and $y$ are semisimple with distinct eigenvalues 
and that, for any eigenvectors $v$ and  $w$ of $x$ and
$y$ respectively, $w^*(v)\ne 0$.

Our infinite order elliptic $\alpha$ clearly lies in $S$. 
Since the closure of $H$ is all of $PSL(2,\bold R)$, $H$
must contain a hyperbolic element. Conjugate $\alpha$ by
this hyperbolic element to get an elliptic $\beta$ with
a different fixed point. $\beta$ is also semi-simple and has the 
same eigenvalues as $\alpha$, so $\beta$ is also in $S$. 
A straightforward  calculation shows that the eigenvectors for
$\alpha$ and $\beta$ satisfy the necessary condition, so 
$(\alpha, \beta)\in\Cal U$. We can now prove:

\proclaim{Proposition \newitem}$H$ contains two small elements $a$ and $b$
that generate a free group.
\endproclaim

\demo{Proof}Since $(\alpha, \beta)$ lies in $\Cal U$, 
there is some integer $m$ such that $\alpha^m$ and $\beta^m$
generate a free group. Furthermore since $\alpha^m$ and $\beta^m$
are infinite order elliptics, there exist integers $k$ and $l$ such that
$\alpha^{km}$ and $\beta^{lm}$ are $\epsilon_0$-close to the identity.

Since a subgroup of a free group is free, $a=\alpha^{km}$ and 
$b=\beta^{lm}$ are small elements generating a free group.
\qed\enddemo

Our goal now is to use $a$ and $b$ to find a large number of
elements near the identity. Let $W_i$ be sets of words in $a$
and $b$ defined inductively as follows. Let $W_0=\{a,b,a^{-1},b^{-1}\}$.
Let $W_i=\{[x,y]|x,y\in W_{i-1}, x\ne y , x\ne y^{-1}\}$.
Let $W=\bigcup W_i$. Let $C_i$ be the image of $W_i$ in the free
group $\langle a,b\rangle$. Set $c_i=|C_i|$.

\proclaim{Lemma \newitem}Distinct words in $W$ have distinct
images in the group $\langle a,b\rangle$.
\endproclaim
\name\comm

\demo{Proof}The idea is to show that any word in $W$ can be 
canonically reconstructed from the reduced word representing 
the same element in the group $\langle a,b\rangle$. 
In order to describe the
reconstruction, it is helpful to describe a step-by-step reduction
of a word in $W$. 

This paragraph explains how, given a word $w\in W_i$, words
$w_{i-1},w_{i-2},\dots,w_0$ are inductively
defined, where $w_k$ is
a partial reduction of $w_{k+1}$. Note the indices decrease
as the induction proceeds.
View $w$ as a word consisting of 4 blocks, each of length
$4^{i-1}$. By construction no
two consecutive blocks are inverses of each other. Set
$w_{i-1}=w$. Given $w_j$, for $j\geq 1$, we construct $w_{j-1}$ by viewing
$w_j$ as being made up of blocks of length $4^{j-1}$ and
canceling blocks as follows. Scan $w_j$ from left to right
until a block followed by its inverse appears. Cancel these
two blocks and then continue scanning from left to right,
starting with the block which followed those just cancelled,
until another block is followed by its inverse. Cancel these
and continue this process until the end of the word is reached.
Call the resulting word $w_{j-1}$. 

Conceivably $w_{j-1}$ contains
consecutive blocks of length $4^{j-1}$ which are inverses of
each other. The next step is to show that there are not any such
pairs. We proceed by induction with the indices decreasing. 
As noted above $w_{i-1}$ doesn't contain any consecutive blocks
of length $4^{i-1}$ which are inverses of each other. 
By induction assume that $w_j$ contains no consecutive
blocks  of length $4^j$ which are inverses of each other. We show that
$w_{j-1}$ contains no consecutive blocks of length $4^{j-1}$
which are inverses of each other. If there were consecutive
blocks of length $4^{j-1}$ which were inverses of each other,
they would have to have come from a sequence in
$w_j$ of the form $x_ix_jx_j^{-1}x_i^{-1}$ where $x_i$ and
$x_j$ are in $W_{j-1}$.  Since $w_j$
was gotten by canceling only blocks of length $4^i$ or longer
this sequence must be part of a sequence of the form
$x_i^{-1}x_j^{-1}x_ix_jx_j^{-1}x_i^{-1}x_jx_i$, but such a
sequence cannot occur since by induction we are assuming that
$w_j$ contains no two  consecutive blocks which are inverses
of each other. Thus $w_{j-1}$ contains no consecutive blocks
of length $4^{j-1}$ which are inverses of each other. In
particular this argument shows that $w_0$ is a reduced word,
since it can contain no consecutive blocks of length $4^0$
which are inverses of each other.

Let $R_0$ be the set of all words gotten by completely reducing some word in 
$W$. Given a word $w_0$ in $R_0$, we wish to canonically 
reconstruct the element of $W$ it came from.
A word $w_0\in R_0$ is either in $W_0$ or it came from
reducing a word $v_1$ with the property that it is made up
of blocks of length 4, where each block is an element of $W_1$,
and no consecutive blocks are inverses of each other. We will
say a word $v$ has property $(\ast)$ if it reduces to $w_0$ and
is made up of blocks of length 4, where each block is an
element of $W_1$, and no consecutive blocks are inverses 
of each other.

This paragraph is devoted to
showing that there is only one possible word satisfying $(\ast)$. First we must show
that any possible word satisfying $(\ast)$ 
must start with the same three letters as
$v_1$ does. The word $v_1$ must start with some sequence of
the form $x_1x_2x_1^{-1}x_2^{-1}$ for some $x_i\in W_0$. If the second block
starts with $x_2$ then when reducing to $w_0$ the $x_2^{-1}x_2$
cancel but no other cancellation takes place between these two blocks
since they aren't allowed to be inverses of each other. Similarly
the last and first letter of consecutive blocks may cancel but
nothing more, so in particular nothing ever cancels with the
first three letters of $v_1$. Thus $w_0$ and $v_1$ must start
with the same three letters. Moreover any word satisfying $(*)$
must start with the same four letters
as $v_1$. If $w_0$ and $v_1$ begin with the same four letters
then let $w_0'=w_0$. Otherwise let $w_0'$ be the word gotten
by inserting the fourth letter of $v_1$ followed by its inverse
into $w_0$ between its third and fourth letter. We have just
shown that the expansion $w_0'$ of $w_0$ is independent of the 
choice of $v_1$. The same reasoning shows that any two elements
satisfying $(*)$ and beginning with the same $4n$ letters which reduce to $w_0$
must agree on the first $4(n+1)$ letters. So by induction
there is only one $w_1$ satisfying $(*)$ which reduces to $w_0$.

Now say $w_i$ has been reconstructed. Then $w_i$ is either
in $W_i$ or it comes from reducing a word $w_{i+1}$ made up
of blocks of commutators of elements of $W_i$ such that no two 
consecutive blocks are inverses of each other. By the same
reasoning we used in the case $i=0$, only one such word exists.
Thus by induction we may canonically reconstruct $w\in W$ from $w_0$.
Thus no two distinct words in $W$ can represent the same element 
of $\langle a,b\rangle$.

\qed\enddemo

We want to count
the elements in $W$ and see how the number grows with the word length.
Note $c_i=c_{i-1}(c_{i-1}-2)\geq\frac{1}{2}c_{i-1}^2.$

\proclaim{Observation \newitem}$c_i\geq2^{2^i}.$
\endproclaim

\demo{Proof}We will prove that $c_i\geq2^{2^i+1}$, from which the
observation follows. For $n=1$, we have $c_1=8=2^{2^1+1}$. Assume
$c_{i-1}\geq2^{2^{i-1}+1}$. Then,
	$$\align
	c_i&\geq\frac{1}{2}(2^{2^{i-1}+1})^2\\
		&=\frac{1}{2}(2^{2^i+2})\\
		&=2^{2^i+1}.
	\endalign$$
\qed\enddemo

\proclaim{Theorem \newitem}For any finitely generated 
subgroup $H$ of $PSL(2,\bold R)$ which is not discrete or elementary,
$|H_{\epsilon_0}^n|$ grows faster than
$f(n)=2^\frac{\sqrt n}{4}$. 
\endproclaim
\name\fastgrowth

\demo{Proof}Elements in $C_i$ have word length no more
than $4^i$. So $|H^{4^i}_{\epsilon_0}|\geq 2^{2^i}$. Given
$n$, let $j$ be such that $4^j\leq n <4^{j+2}$.
	$$f(n)=2^\frac{\sqrt n}{4}<2^{2^j}\leq
|H^{4^j}_{\epsilon_0}|\leq |H^n_{\epsilon_0}|.$$
\qed\enddemo

\proclaim{Theorem \newitem}G must be discrete.
\endproclaim

\demo{Proof}Apply Theorem \fastgrowth\ together with Theorem \lingrowth\
taking $\epsilon=\epsilon_0$. 
\qed\enddemo

\head\newhead. The kernel of $\Phi$ is quasi-isometric to $\bold R$
\endhead

To complete our proof we need only show that the kernel of
$\Phi:\Gamma\to G$ is infinite and therefore, by Lemma \finorcqi, 
quasi-isometric to $\bold R$.

\proclaim{Lemma \newitem}The kernel of $\Phi$ is infinite.
\endproclaim

\demo{Proof}Since $G$ is discrete there are only finitely many elements
$g_u$ that move $x_0$ less than any given bounded amount. For any
$u\in\Gamma$ such that $\phi_u$ moves  $z_0$ less than or equal to
$M_1$ horizontally, the corresponding $g_u$ can move
$x_0$ no more than $M_1+L$. But there are infinitely
many such $u$, since 
any point in $\bold H^2\times\bold R$ is within $M_1$
of some orbit point of $z_0$.  So some element $g\in G$ must be $g_u$ for
infinitely many $u$, say $u_1,u_2,u_3,\dots$. But then all of
the $g_{u^{-1}u_i}$ must be the identity, hence $\ker\Phi$
is infinite.
\qed\enddemo

Thus we have shown:

\proclaim{Theorem \newitem}Any group $\Gamma$ 
quasi-isometric to $\h2r$ is an extension of a finite extension of a
cocompact Fuchsian
group by a virtually infinite cyclic group. In other words there is an
exact sequence
	$$0 @>>> A @>>> \Gamma @>>> G @>>> 0$$
where $A$ is virtually infinite cyclic and $G$ is a finite extension
of a cocompact Fuchsian group.
\endproclaim

\newpage 
\vskip 20pt
\centerline{REFERENCES}
\vskip 20pt

1. A. Beurling and L. Ahlfors, ``Boundary correspondence 
of quasiconformal mappings'', {\it Acta Math.}  96 (1956), 125-142.

2. Y. Carri\`ere, ``Feuilletages riemanniens \`a croissance
polyn\^omiale'', {\it Comment. Math. Helv.}  63 (1988), 1--20.

3. Y. Carri\`ere and \'E. Ghys, ``Relations d'\'equivalence 
moyennable sur les groupes de Lie'', {\it C. R. Acad. Sci. Paris I Math.}
300 (1984), 677--680.

6. A. Casson and D. Jungreis, ``Convergence groups and Seifert
fibered 3-manifolds'', {\it Invent. Math} 118(3) (1994), 441--456.

4. S. S. Chen and L. Greenberg, ``Hyperbolic spaces'' in {\it 
Contributions to analysis}, {\it Acad. Press}  51 (1974), 49--87.

5. R. Chow, ``Groups quasi-isometric to complex
hyperbolic space'', {\it Trans. Amer. Math. Soc.}  348 (1996), 1757--1769.
 
9. V. Efremovitch and E Tichonirova, ``Equimorphisms of hyperbolic
spaces'', Izv. Akad. Nauk CCCP 28 (1964), 1139 -- 1144.

7. A. Eskin, ``Quasi-isometric rigidity of nonuniform lattices in
lattices in higher rank symmetric spaces'', {\it Journal Amer. Math. Soc.}
 10 (1997), 48--80.

8. A. Eskin and B. Farb, ``Quasi-flats and rigidity in higher rank
symmetric spaces'', {\it Journal Amer. Math. Soc.} 10 (3) (1997),
653--692.

10. B. Farb, ``The quasi-isometry classification of lattices in
semisimple Lie groups'', {\it Math. Res. Letters}  4(5) (1997), 705--718.

11. B. Farb and L. Mosher, ``Quasi-isometric rigidity for the solvable 
Baumslag-Solitar groups.'', {\it Inventiones Math.}  131(2)
(1998), 419--451.

12. B. Farb and L. Mosher, ``Quasi-isometric rigidity for the solvable 
Baumslag-Solitar groups, II.'', {\it Inventiones Math.}  137(3)
(1999), 273--296. 

13. B. Farb and R. Schwartz, ``The large-scale geometry of Hilbert
modular groups'', {\it J. Diff. Geom.} 44(3) (1996), 435--478.
 
14. D. Gabai, ``Convergence groups are Fuchsian groups'', {\it Ann. of
Math.} 136 (1992), 447--510.

15. \'E. Ghys, ``Les groupes hyperboliques'', {\it S\'eminaire Bourbaki}
722 (1990), 1--29.

16. \'E. Ghys and P. de la Harpe, eds., {\it Sur les groupes hyperboliques
d'apr\`es Mikhael Gromov}, Prog. in Math  83,  Birkh\"auser (1990).

17. M. Gromov, ``Groups of polynomial growth and expanding maps'',
{\it Inst. Hautes \'Etudes Sci. Publ. Math.}  53 (1981), 53--73.
 
18. M. Gromov, ``Infinite groups as geometric objects'',
{\it Proceedings of the International Congress of Mathematicians, Warsaw 1983},
385-392 (1984).

19. M. Gromov, ``Hyperbolic groups'' in {\it Essays in
group theory}, ed. S. M. Gersten, {\it Math. Sci. Res. Inst.} 
 8 (Springer, 1987), 75--263.
 
20. M. Gromov, `` Asymptotic invariants of infinite groups'',
{\it Geometric Group Theory, Vol. 2} London Math. Soc. Lecture Notes  182,
Cambridge Univ. Press, (1993).

21. A. Hinkkanen, ``The structure of certain quasisymmetric groups'',
{\it Mem. Amer. Math. Soc.}  422 (1990), 1--87.

22. M. Kapovich and B. Leeb, ``Quasi-isometries preserve the geometric 
decomposition of Haken manifolds'', {\it Invent. Math.}
128 (1997), 393--416.

23. B. Kleiner and B. Leeb, ``Rigidity of quasi-isometries for
symmetric spaces and Euclidean buildings'', {\it Pub. IHES}
86 (1997), 115--197.

24. B. Kleiner and B. Leeb, ``Groups quasi-isometric to
symmetric spaces'', preprint (1996).

25. G. Mess, ``The Seifert confecture and groups which are coarse
quasi-isometric to planes'', preprint.

26. G. D. Mostow, {\it Strong rigidity of locally symmetric spaces}, Princeton
University Press (1972).

27. P. Pansu, ``M\'etriques de Carnot-Carathe\'odory et
quasi-isom\'etries des espaces de rang 1'', {\it Ann. of Math.} 129
(1989), 1--60.
 
28. M. S. Raghunathan, {\it Discrete subgroups of Lie groups}, Springer-
Verlag (1972).
 
29. R. Schwartz, ``The quasi-isometry classification of rank one
lattices'', {\it IHES Sci. Publ. Math.},  82 (1996).

30. R. Schwartz, ``Quasi-isometric rigidity and diophantine
approximation'', {\it Acta Math.}  177(1) (1996), 75--112.

31. P. Scott, ``The geometries of three-manifolds'', {\it Bull. London
Math. Soc.}  15 (1983), 401--487.

32. J. Stallings, ``On torsion-free groups with infinitely
many ends'', {\it Ann. of Math.} bf 88 (1968), 312--334.

33. W. P. Thurston, {\it Three-dimensional geometry and topology, Vol 1},
Princeton Univ. Press, (1997).

34. J. Tits, ``Free subgroups of linear groups'', {\it J. Algebra}
20 (1972), 250--270.

35. P. Tukia, ``On quasiconformal Groups'', {\it J. Analyse Math.} 
46 (1986), 318--346.
  
36. P. Tukia, ``Homeomorphic conjugates of Fuchsian groups'', {\it J.
Reine Angew. Math.}  391 (1988), 1--54. 

37. H. Zassenhaus, ``Beweis eines Satzes \"uber diskrete Gruppen'',
{\it Abh. Math. Sem. Hansisch Univ.}  12 (1938), 289--312.

\enddocument

\bye